\documentclass[11pt,a4paper]{article}
\usepackage{a4wide,amsmath,amsthm,epsfig,graphicx}
\usepackage{amsmath,amsthm,amssymb}
\usepackage{amsfonts}
\usepackage{graphicx}
\title{\bf{Implementing Quasi-Monte Carlo Simulations with Linear Transformations}}
\author{Piergiacomo Sabino\\
Dipartimento di Matematica\\
Universit\`{a} degli Studi di Bari\\
sabino@dm.uniba.it\\
Rapporto 41/07}
\date{}
\begin{document}
\maketitle \thispagestyle{empty}
\abstract \noindent Pricing exotic multi-asset path-dependent
options requires extensive Monte Carlo simulations. In the recent
years the interest to the Quasi-monte Carlo technique has been
renewed and several results have been proposed in order to improve
its efficiency with the notion of effective dimension. To this aim,
Imai and Tan introduced a general variance reduction technique in
order to minimize the nominal dimension of the Monte Carlo method.
Taking into account these advantages, we investigate this approach
in detail in order to make it faster from the computational point of
view. Indeed, we realize the linear transformation decomposition
relying on a fast ad hoc QR decomposition that considerably reduces
the computational burden. This setting makes the linear
transformation method even more convenient from the computational
point of view. We implement a high-dimensional (2500) Quasi-Monte
Carlo simulation combined with the linear transformation in order to
price Asian basket options with same set of parameters published by
Imai and Tan. For the simulation of the high-dimensional random
sample, we use a $50$-dimensional scrambled Sobol´ sequence for the
first $50$ components, determined by the linear transformation
method, and pad the remaining ones out by the Latin Hypercube
Sampling. The aim of this numerical setting is to investigate the
accuracy of the estimation by giving a higher convergence rate only
to those components selected by the linear transformation technique.
We launch our simulation experiment also using the standard Cholesky
and the principal component decomposition methods with pseudo-random
and Latin Hypercube sampling generators. Finally, we compare our
results and computational times, with those presented in Imai and
Tan \cite{IT2007}.
\endabstract
%
\section{Introduction}

The Monte Carlo method (MC) is a computational intensive technique
whose purpose is to estimate integrals numerically. It is
characterized by a rate of convergence of order $O(1/\sqrt{n})$,
where $n$ is the number of simulations, and it is independent of the
problem dimension $d$. This last feature makes the MC method
appealing and applicable to several financial high-dimensional
situations such as options pricing. Furthermore, the estimation
error (RMSE), that can be easily computed statistically, depends
only on the convergence rate and on an intrinsic constant.

Based on probabilistic considerations, standard reduction techniques
can only reduce the constant but cannot improve the convergence
rate.

In contrast, Quasi-Monte Carlo methods (QMC) aim to enhance the
convergence rate by means of low-discrepancy sequences. These
sequences provide better stratification and a convergence rate of
order $O\left(\frac{ln^dn}{n}\right)$ (see Niederreiter
\cite{Ni1992}). The rate is faster than the previous one but depends
on the problem dimensions. These sequences are purely deterministic,
meaning that the estimation error cannot be estimated statistically.
In the Randomized Quasi-Monte Carlo (RQMC) method some randomness is
introduced in the low-discrepancy sequences while preserving their
better convergence rate. This technique is called scrambling.

Several numerical investigations conclude that QMC and RQMC
simulations do not give substantial advantage for $d>10/20$.

Some approaches have been proposed in order to extend the QMC
superiority to high-dimensional estimations. Caflisch \emph{et al}
\cite{CMO1997} address the problem using the analysis of variance
(ANOVA) of the integrand function and defining two notions of
effective dimension: the effective dimension in truncation and
superposition sense. Briefly, the truncation dimension reflects
that, for some integrand functions, only a small number of inputs
really matter. The definition of effective dimension in
superposition sense takes into account that for some integrands the
inputs might influence the outcome through their joint action within
small groups.

Imai and Tan \cite{IT2002} proposed a general linear transformation
construction (LT) to reduce the effective dimension of the problem
in superposition sense, focusing on the particular payoff function.
The authors show that this approach offers a considerable advantage
with respect to the principal component analysis (PCA) in terms of
accuracy and versatility.

Moreover, their simulation procedure relies on the complete Latin
Supercube Sampling generation (see Owen \cite{ow1998B} for more on
this topic) in order to generate a high-dimensional low discrepancy
sequence with good properties.

Here we investigate the accuracy of the LT method in detail and
implement the construction fast by an efficient QR decomposition. We
run our simulation procedure with the same set of parameters as in
Imai and Tan \cite{IT2007} and are thus able to directly compare the
respective results.

We will demonstrate, that our implementation makes the LT
considerably faster and maintain its versatility.

We test the efficiency of the LT construction by launching a MC
simulation in a more extreme setting. We use scrambled
low-discrepancy sequencies only to those components the LT considers
as optimal, while simulating the others with the Latin Hypercube
Sampling (LHS) that has lower convergence rate.

The LHS is supposed to give good accuracy when the target function
is a sum of one-dimensional ones. If the LT accomplishes this task
optimally it would give a good improvement in this setting too. Our
experiment is intended to test if the LT gives the same results as
in Imai and Tan \cite{IT2007}, in terms of RMSEs, in this partial
RQMC setting. This means, that if the LT with RQMC provides
considerable advantage with respect to the pure LHS generation it
reduces the effective dimension in superposition sense really
optimally.

As a comparison, we launch the MC simulation using a standard
pseudo-random generator and build the random path with standard
Cholesky and PCA decompositions too.

The paper is organized as follows. Section 2 describes the financial
setting and formulates the Asian basket option pricing problem as an
integral explicitly. Section 3 introduces the MC and the QMC methods
and the notion of effective dimensions of the problem. Section 4
describes the LT construction introduced by Imai and Tan and how it
applies to several financial situations. Section 5 presents the main
steps of our MC simulation. Section 6 illustrates the numerical
results we obtain and discuss the efficiency of the LT and its fast
implementation. Section 7 concludes the paper and the Appendix
describes the \emph{ad hoc} QR decomposition used.

%
\section{Problem Statement}
We consider the problem of estimating the fair price of a contract
in a standard financial market $\mathfrak{M}$ in a Black-Scholes
framework, with  a constant risk-free rate $r$ and time-dependent
volatilities. There are $M + 1$ assets in the market, one risk free
asset and $M$ risky assets. The price processes of the assets in
this market are driven by a set of stochastic differential
equations.

Suppose we have already applied the Girsanov
theorem and found the (unique) risk-neutral probability, the model
for the risky assets is the so called multi-dimensional geometric
brownian motion:
\begin{equation}
S_{0}(t)=e^{rt}\label{1.1}
\end{equation}
\begin{equation}
dS_{i}\left( t\right) =rS_{i}\left( t\right) dt+\sigma _{i}\left(
t\right)S_{i}\left( t\right)\,dW_{i}\left( t\right) ,\qquad
i=1,\dots ,M.  \label{1.2}
\end{equation}%
\noindent Here $S_{i}\left( t\right) $ denotes the $i$-th asset
price at time $t$, $\sigma _{i}\left( t\right)$ represents the
instantaneous time-dependent volatility of the $i$-th asset return,
$r$ is the continuously compounded risk-free  interest rate, and $%
\mathbf{W}\left( t\right) =\left( W_{1}\left( t\right) ,\dots
,W_{M}\left( t\right) \right) $ is an $M$-dimensional Brownian
motion. Time $t$ can vary in $\mathbb{R}_{+}^*$, that is, we can
consider any maturity $T\in\mathbb{R}_{+}^*$ for all financial
contracts.

The multi-dimensional brownian motion $\mathbf{W}\left( t\right)$ is
a martingale, each component
is a martingale, and satisfies the following properties:%
\begin{equation*}
\mathbb{E}\left[ W_{i}\left( t\right) \right] =0,\qquad i=1,\dots
,M.
\end{equation*}%
\begin{equation*}
\left[W_{i},W_{k}\right]\left(t\right)=\rho _{ik}t,\qquad i,k =
1,\dots ,M.
\end{equation*}%
\noindent where $[\centerdot,\centerdot](t)$ represents the
quadratic variation up to time $t$ and $\rho _{ik}$ the constant
instantaneous correlation between $W_{i}$ and $W_{k}$.

Applying the risk-neutral pricing formula, the value at time $t$ of
any European $T$-maturing derivative contract is:
\begin{equation}
V(t) =
exp\left(r(T-t)\right)\mathbb{E}\left[\phi(T)\right|\mathcal{F}_t]\label{1.3}\text{.}
\end{equation}
\noindent $\mathbb{E}$ denotes the expectation under the risk
neutral probability measure and $\phi(T)$ is a generic
$\mathcal{F}_T$ measurable function, with
$\mathcal{F}_T=\sigma\{0<t\leq T;\textbf{W}(t)\}$, that determines
the payoff of the contract. Although not explicitly written, the
function $\phi(T)$ depends on the entire multi-dimensional brownian
path up to time $T$.

We will restrict our analysis to Asian options that  are exotic
derivative contracts that can be written both on a single security
and on a basket of underlying securities. Hereafter we will consider
European-style Asian options whose underlying securities coincide
with the $M+1$ assets on the market. This is the most general case
we can tackle in the market $\mathfrak{M}$ because it is complete in
the sense that we can hedge any financial instrument by finding a
portfolio that is a combination of this $M+1$ assets.
\subsection {Asian Options Payoff}
The theoretical price for a discretely monitored Asian option is:
%
%
\begin{equation}
a_i\left( t\right) =exp\left(r(T-t)\right)\mathbb{E}\left[\left(
\frac{\sum_{j=1}^{N}S_{i}\left( t_{j}\right)}{N}
-K\right)^+\bigg|\mathcal{F}_t\right] \quad\text{Option on a
Single Asset}\label{1.1.3}
\end{equation}
\begin{equation}
a\left( t\right) =exp\left(r(T-t)\right)\mathbb{E}\left[\left(
\sum_{i=1}^{M}\sum_{j=1}^{N}w_{ij}\,S_{i}\left( t_{j}\right)
-K\right)^+\bigg|\mathcal{F}_t\right]\quad\text{Option on a
Basket} \label{1.1.4}
\end{equation}%
\noindent where $t_1<t_2\dots<t_N=T$ and the coefficients $w_{ij}$
satisfy $\sum_{i,j}w_{ij}=1$.

European options with payoff functions $\left( \ref{1.1.3}\right)$
and $\left( \ref{1.1.4}\right)$ are called arithmetic weighted
average options or simply arithmetic Asian options. When $M>0$ and
$N=1$ the payoff only depends on the terminal price of the basket of
$M$ underlying assets and the option is known as basket option.

\subsection {Problem Formulation as an Integral}
The model $\mathfrak{M}$, presented in the first section, consists
of the risk-free money market account and $M$ assets driven by $M$
geometric brownian motion described by equation (\ref{1.2}) whose
solution is:
\begin{equation}
S_{i}\left( t\right) =S_{i}\left( 0\right) exp\left[
\int_0^t\left( r- \frac{\sigma _i^2\left( s\right)}{2}\right)ds
+\int_0^t\sigma _i\left( s\right)dW_{i}\left( s\right) \right]
,i=1,...,M.\label{1.2.1}
\end{equation}
\noindent The quantity $\int_0^t\frac{\sigma _i^2\left(
s\right)}{2}ds$ is the total volatility for the $i$-th asset. The
solution (\ref{1.2.1}) is a multi-dimensional geometric brownian
motion, written GBM$\left(r,\int_0^t\frac{\sigma _i^2\left(
s\right)}{2}ds\right)$.

Under the assumption of constant volatility the solution is still a
multi-dimensional geometric brownian motion with the following form:
\begin{equation}
S_{i}\left( t\right) =S_{i}\left( 0\right) exp\left[ \left( r-
\frac{\sigma _i^2}{2}\right)t +\sigma _iW_{i}\left( t\right)
\right] ,i=1,...,M.\label{1.2.2}
\end{equation}
\noindent In compacted notation the solution (\ref{1.2.2}) is
GBM$\left(r,\frac{\sigma _i^2}{2}t\right)$.

Pricing Asian option requires to monitor the solutions (\ref{1.2.1})
and (\ref{1.2.2}) at a finite set of points in time
$\{t_1,\dots,t_N\}$. This sampling procedure yields to the following
expressions for time-depending and constant volatilities:
\begin{equation}
S_i(t_j) = S_i(0)exp\bigg[\int_0^{t_j}\left(r -
\frac{\sigma^2_i(t)}{2}\right)dt + Z_i(t_j)\bigg]\label{1.2.3}
\end{equation}
\begin{equation}
S_i(t_j) = S_i(0)exp\bigg[\left(r - \frac{\sigma^2_i}{2}\right)t_{j}
+ Z_i(t_j)\bigg]\label{1.2.4}
\end{equation}
\noindent where the components of the vector $\left(Z_1(t_1),\dots
Z_1(t_N),Z_2(t_1),\dots,Z_M(t_N)\right)$ are  $M\times N$ normal
random variables with zero mean vector and the following
covariance matrix:
\begin{equation}\label{Cov1}
\Sigma_{MN} = \left( \begin{array}{cccc} \Sigma(t_1) & \Sigma(t_1)
& \ldots & \Sigma(t_1)
\\ \Sigma(t_{1}) & \Sigma(t_2)& \ldots & \Sigma(t_2)
\\ \vdots & \vdots & \ddots & \vdots
\\\Sigma(t_{1}) & \Sigma(t_2) & \ldots & \Sigma(t_{N})
\end{array} \right)\text{for time-dependent volatilities,}
\end{equation}
\noindent or
\begin{equation}
\Sigma_{MN} = \left( \begin{array}{cccc} t_{1}\Sigma & t_{1}\Sigma
& \ldots & t_{1}\Sigma
\\ t_{1}\Sigma & t_2\Sigma & \ldots & t_{2}\Sigma
\\ \vdots & \vdots & \ddots & \vdots
\\t_{1}\Sigma & t_2\Sigma & \ldots & t_{N}\Sigma
\end{array} \right)\label{Cov2}\text{for constant volatilities.}
\end{equation}
Each element depends on four indexes:
\begin{equation}
\Big(\big(\Sigma_{MN}\big)_{ik}\Big)_{lm} = \int_0^{t_l\wedge t_m}
\sigma_i(t) \sigma_k(t)\rho_{ik}dt\label{par}
\end{equation}
with $i,k=1,\dots,M$ and $l,m=1,\dots,N$.

The payoff at maturity $T$ of the arithmetic average Asian option
is then:
\begin{equation}
p_a(T) = \left(g(\mathbf{Z})-K\right)^+\label{1.2.5}
\end{equation}
where
\begin{equation}
g(\mathbf{Z})=\sum_{k=1}^{M\times N} exp\left (\mu_k + Z_k
\right)\label{1.2.6}
\end{equation}
and
\begin{equation}
\mu_k = \ln(w_{k_1k_2}S_{k_1}(0)) +
\bigg(r-\frac{\sigma_{k_1}^2}{2} \bigg)t_{k_2}\label{1.2.7}
\end{equation}
for constant volatilities or
\begin{equation}
\mu_k = \ln(w_{k_1k_2}S_{k_1}(0)) +
rt_{k_2}-\frac{\int_0^{t_{k_2}}
\sigma_{k_1}^2(t)dt}{2}\label{1.2.8}
\end{equation}
for time-dependent volatilities. The indexes $k_1$ and $k_2$ are
respectively $k_1=(k-1)modM, k_2 = [(k-1)/M]+1$, where $mod$ denotes
the modulus and $[\centerdot]$  the greatest integer less than or
equal to $x$.

The calculation of the price $a(t)$ in equation ($\ref{1.1.4}$)
can be formulated as an integral on $[0,1]^{NM}$ in the following
way (see Dahl and Benth \cite{DB2001} and \cite{DB2002}):
\begin{equation}
    a\left( t\right) =exp\left(r(T-t)\right)
    \int_{[0,1]^{NM}}\left(g(\mathbf{u})-K\right)^+F^{-1}_{\mathbf{Z}}(\mathbf{u})\mathbf{du}\label{2.2.10}
\end{equation}

\section{Problem Dimension}

The main purpose of the standard MC method is to numerically
estimate the following integral:
\begin{equation}
I = \int_{[0,1]^d} f(\mathbf{x})\,\mathbf{dx}
\label{3.1.5}\text{.}
\end{equation}
\noindent $I$ can be seen as $\mathbb{E}\left[f(U)\right]$, the
expected value of a function $f(\centerdot)$ of the random vector
$\mathbf{U}$ that is uniformly distributed in hypercube $[0,1]^d$.

MC methods simply estimate $I$ by drawing a sample of $n$
independent replicates $U_1\dots,U_n$ of $\mathbf{U}$ and then
computing the arithmetic average:
\begin{equation}
\widehat{I} = \widehat{I}_n = \frac{1}{n} \sum_{i=1}^n
f(U_i).\label{3.1.6}
\end{equation}

The Law of Large Numbers ensures that $\widehat{I}_n$ converges to
$I$ in probability almost surely and the Central Limit Theorem
states that $I - \widehat{I}_n$ converges in distribution to a
normal with mean $0$ and standard deviation $\sigma/\sqrt{n}$ with
$\sigma=\sqrt{\int_0^1\left(f(\mathbf{x})-I\right)^2\mathbf{dx}}$.
The convergence rate is than $O(1/n)$ for all dimensions $d$. The
parameter $\sigma$ is generally unknown in a setting in which $I$ is
unknown, but it can be estimated using the sampled standard
deviation or root mean square error (RMSE):
\begin{equation}
 RMSE = \sqrt{\frac{1}{n-1}\sum_{i=1}^n\left( f(U_i) -
 \widehat{I}_n\right)^2}.\label{3.1.7}
\end{equation}
When the nominal dimension $d$ of the problem of estimating the
integral $ (\ref{3.1.5})$ is one, there are standard numerical
techniques that give a good accuracy when $f$ is smooth.
Considerable problems arise when $d$ is high.

We aim to estimate the fair value of the Asian option of equation
(\ref{2.2.10}) with an high-dimensional Quasi MC simulation as
formulated in (\ref{3.1.5}).

QMC method relies on the construction of deterministic sequences,
also known as low-discrepancy sequences, that cover the hypercube
$[0,1)^d$ uniformly. We define the quantity $D_{n}^{\ast
}=D_{n}^{\ast }\left( P_{1},\dots ,P_{n}\right) $ as the star
discrepancy. It is a measure of the uniformity of the sequence
$\left\{ P_{n}\right\} _{n\in \mathbb{N}^{\ast }}\in \left[
0,1\right) ^{d}$ and it must be stressed that it is an analytical
quantity and not a statistical one.

A sequence $\left\{ P_{n}\right\} _{n\in \mathbb{N}^{\ast }}$ is
called low-discrepancy sequence if:
\begin{equation}
D_{n}^{\ast }\left( P_{1},\dots ,P_{n}\right) =O\left( \frac{\left(
\ln n\right) ^{d}}{n}\right) .  \label{6.2.1}
\end{equation}

The following inequality, attributed to Koksma and Hlawka, provides
an upper bound to the estimation error of the unknown integral with
the QMC method in terms of the star discrepancy:
\begin{equation}
|I-\hat{I}|\leq D_{n}^{\ast }\,\,V_{HK}\left( f\right) .
\label{6.2.2}
\end{equation}%
\noindent $V_{HK}\left( f\right) $ is the variation in the sense of
Hardy and Krause. Consequently, if $f$ has a finite variation and
$n$ is large enough, the QMC approach gives an error smaller than
the error obtained by the crude \emph{MC} method for low dimensions
$d$.

It is well known that QMC methods loose the better accuracy in high
dimension. It is then fundamental to capture the most important (in
statistical sense) components or to reduce the nominal dimension of
the problem by means of ANOVA considerations.

Let $\mathcal{A}=\{1,\dots,d\}$ denote the set of the independent
variables for $f$ on $[0,1]^d$. $f$ could be written into the sum of
orthogonal functions each of them defined in a different subset of
$\mathcal{A}$, that is depending only on the variables in each of
these subsets:
\begin{equation}
f(\mathbf{x}) = \sum_{u\subseteq \mathcal{A}} f_u(\mathbf{x})\label{3.1}
\end{equation}
Now  let $|u|$ denote the cardinality of $u$ and $\sigma^2 =
\int(f(\mathbf{x})-I)^2\,\mathbf{dx}$, $\sigma^{2}_{u} = \int
f_u(\mathbf{x})^2\,\mathbf{dx}$, $\sigma^{2}_{0}=0$, supposing
$\sigma < +\infty$ and $|u|> 0$ it holds:
\begin{equation}
\sigma^2 = \sum_{u\subseteq \mathcal{A}} \sigma^{2}_{u}\label{3.3}
\end{equation}
Equation (\ref{3.3}) partitions the total variance into parts
corresponding to each subset $u\subseteq\mathcal{A}$. The $f_u$
enjoys some nice properties: if $j\in u$ the line integral
$\int_{[0,1]} f_u(\mathbf{x})\,dx_j = 0$ for any $x_k$ with $k\neq
j$, and if $u\neq v$  $\int f_u(x) f_v(x)\,dx = 0$.

Exploiting the ANOVA decomposition, the definition of
effective dimension can be given in the following ways:
\newtheorem{Definition1}{Definition}
\begin{Definition1}
The effective dimension of $f$, in the superposition sense, is the
smallest integer $d_S$ such that $\sum_{0<|u|\leq d_S} \sigma^2_u
\geq p \sigma^2$ .\\ The value $d_S$ depends on the order in
which the input variables are indexed.
\end{Definition1}
\newtheorem{Definition2}[Definition1]{Definition}
\begin{Definition2}
The effective dimension of $f$, in the truncation sense, is the
smallest integer $d_T$ such that $\sum_{u\subseteq\{1,\dots,d_t\}}
\sigma^2_u \geq p \sigma^2$.
\end{Definition2}
\noindent $0<p<1$ is an arbitrary level; the usual choice is $p=99\%$.

The definition of effective dimension in truncation sense reflects that for some integrands,
only a small number of the inputs might really matter. The
 definition of effective dimension in superposition sense takes into account that for some
integrands, the inputs might influence the outcome through their
joint action within small groups. Direct computation leads to: $d_S
\leq d_T \leq d$.

%
\section{Linear Transform Construction}
Imai and Tan \cite{IT2002} proposed a general LT method for path
generation with main purpose to minimize the effective dimension in
truncation sense of a simulation problem.

The LT approach provides the same results as the PCA-based one,
moreover proves to be more accurate and versatile in certain
situations.

Many studies demonstrate that the QMC pricing of certain specific
derivative contracts is not substantially improved by the brownian
bridge construction. This suggests to focus the attention onto the
particular payoff function while even the PCA approach is applicable
only for multi-dimensional normal random variables. In contrast, the
LT generation focuses on the particular payoff function instead of
the  multi-dimensional brownian path.

This method provides the best results for linear combinations of
normal random variables. Imai and Tan \cite{IT2002}, \cite{IT2005}
and \cite{IT2007} investigated the practical improvement of the LT
method by running very high-dimensional simulations for European
options, bonds pricing in different dynamics (see the cited
references for more on this topic).

A $n$-dimensional random vector $\textbf{Y}$ with covariance matrix
$\Sigma_y$ can be characterized starting from a vector of
independent standard normal variables $\epsilon$ by the following
transformation: $y = C\epsilon$, with $C C^T = \Sigma_y$. Imai and
Tan consider the following class of LT as  solution of the previous
general problem:
\begin{equation}\label{4.3.1}
    C^{LT} = C^{Ch} A
\end{equation}
\noindent where $C^{Ch}$ is the Cholesky matrix associated to the
covariance matrix of the normal random vector to be generated and
$A$ is an orthogonal matrix, i.e. $AA^T=I$.

The optimum $C^{LT}$ is obtained by optimally choosing $A$ so that
the effective dimension in the truncation sense of the problem of
interest is minimized.

Maximizing the explanatory variability of a normal vector with
covariance matrix $\Sigma$ consists in finding the optimum
orthogonal matrix $A^*$  by iteratively solving the following
optimization problem:
\begin{equation}\label{4.3.2}
    \max \|\mathbf{C_{\cdot k}^{LT}}\|^2 = \max_{\mathbf{A_{\cdot k}} \in \mathbb{R}^{NM}}
    \sum_{p=1}^{MN} \mathbf{C_{p\cdot}^{Ch}}\mathbf{A_{{\cdot} k}}
\end{equation}
\noindent subject to $\|\mathbf{A_{\cdot k}}\|=1$ and
$\mathbf{A_{\cdot k}}\cdot \mathbf{A_{\cdot i}^*}=0$ for
$i=1,\dots,k-1$ and $k\le n$ ($\mathbf{A_{\cdot i}^*}$ indicates the
columns that have been already calculated)

$\mathbf{C_{\cdot k}^{LT}}$ represents the $k$-th column vector and
$\mathbf{C_{k \cdot}^{LT}}$ the $k$-th row vector of $C^{LT}$; the
same notation holds for all the matrices. Imai and Tan \cite{IT2007}
proves that this procedure achieves the same results, in terms of
explained variability, of the PCA decomposition of $\Sigma_{MN}$.
Indeed:
\begin{equation}\label{4.2.3}
    \max \|\mathbf{C_{\cdot k}^{LT}}\|^2 = (\mathbf{C^{Ch}}\mathbf{A_{\cdot
    k}})^T \mathbf{C^{Ch}}\mathbf{A_{\cdot k}} = \mathbf{A_{\cdot k}^T}\Sigma \mathbf{A_{\cdot k}}
\end{equation}
Hence the optimization problem is similar to seeking the $k$-th
principal component.

Finding the optimal matrix $A$ is equivalent of finding the optimal
QR transformation of $C$ with $CC^T = \Sigma$ where $R = (C^{Ch})^T$
and $Q = A$ in the sense described before.

The PCA decomposition provides the best solution for normal random
vectors with $Q = V^T$ and $R = \Lambda^{1/2}$ with $V$ and
$\Lambda$ the orthogonal matrix of the eigenvectors and $\Lambda$
the diagonal matrix of all the eigenvalues in decreasing order
respectively.

\subsection{Special Cases}
As for linear combinations of normal random variables the LT
approach minimizes the effective dimension in truncation sense. It
is established from standard statistics that a linear combination of
normal random variables is still a normal random variable with mean
and variance that depend on the linear combination. It is than
trivial that an integral problem with the nominal dimension $d$ that
involves a linear combination of $d$ normal random variables has an
effective dimension in superposition sense equal to one. The LT
procedure returns this results in truncation sense as an
optimization procedure.

Let $f(\mathbf{z})$ be a linear combination of $d$ normal random
variables $f(\mathbf{z}) =\sum_{i=1}^d w_iz_i$, with $\mathbf{z}\sim
N(\mathbf{\mathbf{\mu}};\Sigma)$ and constants $w_i$, $i=1,\dots,d$.
If $C$ denotes the generic decomposed matrix of $\Sigma$  then the
above function can be expressed as:
\begin{equation}\label{4.3.4}
    f(\mathbf{\epsilon}) =\sum_{k=1}^d \alpha_k \epsilon_k + \mathbf{\mu}\cdot \mathbf{w}
\end{equation}
\noindent where $\alpha_k= \mathbf{C_{\cdot k}}\cdot \mathbf{w}$ and
$\mathbf{\epsilon}$ is a $d$-dimensional vector of standard and
independent normal random variables. Furthermore the total variance
of $f$ is:
\begin{equation}\label{4.3.5}
    \sigma^2 = \sum_{i=k}^d \alpha_k^2
\end{equation}
\noindent The truncation dimension is the smallest integer $d_T$
that satisfies:
\begin{equation}\label{4.3.6}
    \sum_{i=k}^{d_T} \alpha_k^2 \ge p\sigma^2
\end{equation}
As with the LT approach the optimal $C$ is $C^{LT}=C^{Ch}A$ that
leads to:
\begin{equation}
    \alpha_k = \mathbf{A_{\cdot k}}\cdot\mathbf{B}\quad k=1\dots,d\label{4.3.7}
\end{equation}
\noindent where $\mathbf{B}=(C^{Ch})^T\mathbf{w}$. Consequently,
minimizing the effective dimension in the truncation sense is
equivalent to maximizing the variance contribution due to the first
component $\alpha_1^2$ and obtaining $\mathbf{A_{\cdot 1}}$.
Iterating this procedure and imposing the orthogonality condition we
get the optimal matrix $A$. It can be proven, see Imai Tan
\cite{IT2002} or \cite{IT2007}, that the optimal solution for $k=1$
is:
\begin{equation}\label{4.3.8}
    \mathbf{A_{\cdot 1}^*} = \pm\frac{\mathbf{B}} {\|\mathbf{B}\|}
\end{equation}
\noindent while for $k=1,\dots,d$ the column vectors can be
arbitrary but must satisfy the orthogonality condition.
Substituting this results into equation (\ref{4.3.6}) we are left
with $\alpha_1=\pm \|\mathbf{B}\|$ and $\alpha_k=0$ for
$k=2,\dots,d$. The original function $f$ can be written as:
\begin{equation}\label{4.3.9}
    f(\mathbf{\epsilon})= \mathbf{\mu}\cdot\mathbf{w}\pm\|\mathbf{B}\|\epsilon_1
\end{equation}
\noindent This is the best possible scenario for the dimension
reduction. The LT approach reduces any nominal $d$-dimensional
problem involving a linear combination of normal random variables
into a one-dimensional problem in truncation sense. This means that
the LT method rearranges the linear structure of the function for
the best possible reduction.

Let us now consider the following function:
\begin{equation}\label{4.3.10}
    f(\mathbf{\epsilon}) = exp\left(\mu + \sum_{k=1}^n \alpha_k\epsilon_k
    \right) - K
\end{equation}
\noindent with $\mu$, $\alpha_k$ and $K$ constant.

$f(x)^+$ can be considered the payoff function of a geometric
average Asian option with strike price $K$ and:
\begin{equation}\label{4.3.11}
\begin{array}{cc}
  \mu = \sum_{i=1}^M\sum_{j=1}^Nw_{ij}\left[log
    S_i(0)+\left(r-\frac{\sigma_i^2}{2}t_j\right)\right], & \alpha_k
    = \mathbf{C_{\cdot k}}\cdot \mathbf{w}
\end{array}.
\end{equation}
Such a derivative contract is not traded but nevertheless serves to
understand the computational problem. Indeed, performing the
logarithm $log(f\mathbf(\epsilon)-K)$, we obtain a new function
which is a linear combination of normal variates. Applying the
results of the LT method for the previous example we showed that the
nominal dimension $MN$ of the new problem shrinks to one. Again this
is not surprising because we know that the product of log-normal
variates is still a log-normal variate.

These examples highlight the main differences between the PCA
decomposition and the LT methods. The former returns the best
decomposition of the covariance matrix of a normal random vector in
terms of variability of each component. The latter reduces the
effective dimension of the problem focusing on the particular payoff
function. It provides the best solution for linear combinations of
normal variates.
\subsection{General Case}
General payoff functions for European style options are neither
linear combinations of normal random variables nor they can be
obtained by monotone transformations as for the case of the
geometric average Asian options. To address the problem Imai and Tan
propose to approximate an arbitrary  function $g$, such that $g^+$
is the payoff function of a European derivative contract, with its
first order Taylor expansion:
\begin{equation}\label{4.3.12}
    g(\mathbf{\epsilon}) = g(\mathbf{\hat{\epsilon}}) +
    \sum_{l=1}^n\frac{\partial g}{\partial\epsilon_l}\Big|_{\mathbf{\epsilon}=\mathbf{\hat{\epsilon}}}\Delta\epsilon_l
\end{equation}

The approximated function is linear in the standard normal random
vector $\mathbf{\Delta\epsilon}$ and we can rely on the same results
obtained in the previous subsection. By considering an arbitrary
point of expansion, such as $\mathbf{\hat{\epsilon}=0}$, we can
derive the first column of the optimal orthogonal matrix $A^*$. We
can find the complete matrix by expanding $g$ at different points
and then run the optimization algorithm.

Summarizing the optimization can be formulated as follow:
\begin{equation}\label{4.3.13}
    \max_{\mathbf{A\cdot k}\in \mathbf{R^n}}\left(\frac{\partial g}
    {\partial\epsilon_l}\Big|_{\mathbf{\epsilon}=\mathbf{\hat{\epsilon}}}\right)^2
\end{equation}
subject to $\|\mathbf{A_{\cdot k}}\|=1$ and $\mathbf{A^*_{\cdot
j}}\cdot \mathbf{A_ {\cdot k}}=0, j=1,\dots,k-1, k\le n$.

Although equation (\ref{4.3.8}) provides an easy solution at each
step, the correct procedure requires that $\mathbf{A_{\cdot k}}$
must be orthonormal to all the previous (and future) columns. This
feature can be easily obtained by the Gram-Schmidt
orthonormalization or even better by the QR method that is
numerically stabler. As for the latter we must note that the QR
method might return opposite matrices at different time steps
(cosmetic sign adjustment). This does not affect the problem because
the solution in equation (\ref{4.3.8}) can be either with a positive
and negative sign. Furthermore, we stress that it is not necessary
to run the complete QR method at each step. Indeed, all the columns
already calculated are orthogonal and we should use a "partial" QR
method that considerably reduces the computational burden as it will
be shown in the appendix.

Imai and Tan set
$\hat{\mathbf{\epsilon_1}}=\mathbf{0},\hat{\mathbf{\epsilon_2}}=(1,0,\dots,0),\dots,
\hat{\mathbf{\epsilon_k}}=(1,1,1,\dots,0,\dots,0),\dots,\hat{\mathbf{\epsilon_n}}=(1,,\dots,1,0)$',
the $k$-th point has k-1 leading ones.

The choice is arbitrary and a different set can be used that would
return different optimal orthogonal matrices.

Moreover the computational cost can be reduced by only seeking a
suboptimal matrix with optimal columns up to $k^*<n$. This
approximation is reasonable since in practice only a few components
are of relevance as will be shown in the numerical examples.

\subsection{Asian Options Case}

We consider the function $\bar{g}=g-K$ in equation (\ref{1.2.6}), it
is then, easy to verify that its variance can be expressed as:

\begin{equation}\label{4.4.1}
    \sigma^2\left(\bar{g}\left(\mathbf{\epsilon}\right)\right) =
    \sum_{i=1}^{MN}\sum_{j=1}^{MN}exp\left(\mu_i+\mu_j +(1/2)\sum_{l=1}^{MN}\left(C_{il}^2+C_{jl}^2\right)\right)
    \left[exp\left(\sum_{l=1}^{MN}C_{il}C_{jl}\right)-1\right]
\end{equation}

Due to the tractability of the function above, Imai and Tan provide
some implementations of the LT construction. We only show two of
them.

The variance contribution for the first $p$ dimensions can be
defined as:
\begin{equation}\label{4.4.2}
    \sigma_p^2\left(\bar{g}\left(\mathbf{\epsilon}\right)\right) =
     \sum_{i=1}^{MN}\sum_{j=1}^{MN}exp\left(\mu_i+\mu_j +(1/2)\sum_{l=1}^{p}\left(C_{il}^2+C_{jl}^2\right)\right)
    \left[exp\left(\sum_{l=1}^{p}C_{il}C_{jl}\right)-1\right]
\end{equation}
Working with algebra and approximating the exponential in the square
bracket up to the first order, we can obtain the first formulation
for the optimal matrix $A$:
\begin{equation}\label{4.4.3}
    \max_{\mathbf{A_p\cdot}\in\mathbf{R^{MN}}}=
    \sum_{i=1}^{MN}\sum_{j=1}^{MN}exp\left(\mu_i+\mu_j +(1/2)\sum_{l=1}^{p}\left((C_{il}^*)^2+(C_{jl}^*)^2\right)\right)
    C_{il}C_{jl}
    \end{equation}
\noindent subject to $\|\mathbf{A_{\cdot p}}\|=1$ and
$\mathbf{A_{\cdot j}}\cdot\mathbf{A_{\cdot p}^*} = 0$ for $j<k$.

The second formulation consists in applying the general approach by
expanding the function of equation (\ref{1.2.6}) up to the first
order:
\begin{equation}\label{4.4.4}
    g(\mathbf{\epsilon}) = g(\mathbf{\hat{\epsilon}}) +
    \sum_{l=1}^{NM}\left(\sum_{i=1}^{NM}exp\left(\mu_i+\sum_{k=1}^{NM}C_{ik}\hat{\epsilon_k}\right)C_{il}\right)\Delta\epsilon_l
\end{equation}
\noindent We start the optimization procedure by finding the first
column of the optimal matrix $A$:
\begin{equation}\label{4.4.5}
    g(\mathbf{\epsilon}) = g(\mathbf{0}) +
    \sum_{l=1}^{NM}\left(\sum_{i=1}^{NM}exp\left(\mu_i\right)C_{il}\right)\Delta\epsilon_l
\end{equation}
\noindent we set
$\alpha_i=\left(\sum_{i=1}^{NM}exp\left(\mu_i\right)C_{il}\right)=
\sum_{m=1}^{NM}\left(\sum_{i=1}^{NM}exp\left(\mu_i\right)C_{im}^{Ch}\right)A_{ml}$
in order to formulate equation (\ref{4.4.2}) as (\ref{4.3.4}). Set
$\mathbf{d^{(1)}} = (e^{\mu_1},\dots,e^{\mu_{MN}})^T$ and
$\mathbf{B^{(1)}}= \mathbf{(d^{(1)})^T}C^{Ch}$ we know from the
linear combination case that $\mathbf{A_{\cdot
1}^*}=\pm\frac{\mathbf{B^{(1)}}}{\|\mathbf{B^{(1)}}\|}$.

 The $p$-th optimal column can be found considering
the $p$-th starting point of the Imai and Tan's strategy. This
results in:
\begin{equation}\label{4.4.6}
    g(\mathbf{\epsilon}) = g(\mathbf{\hat{\epsilon_p}}) +
    \sum_{l=1}^{NM}\left(\sum_{i=1}^{NM}exp\left(\mu_i+\sum_{k=1}^{p-1}C_{ik}^*\right)C_{il}\right)\Delta\epsilon_l
\end{equation}
\noindent where $C_{ik}^*$, $k<p$ have been already found at the
$p-1$ previous steps and $\mathbf{A_{\cdot p}}$ must be orthogonal
to all the other columns. As for the first step we define
$\mathbf{d^{(p)}} =
\left(exp\left(\mu_1+\sum_{k=1}^{p-1}C_{1k}^*\right),\dots,exp\left(\mu_{MN}+\sum_{k=1}^{p-1}C_{MNk}^*\right)^T\right)$
and $\mathbf{B^{(p)}}= \mathbf{(d^{(p)})^T}C^{Ch}$ the solution is
$\mathbf{A_{\cdot
p}^*}=\pm\frac{\mathbf{B^{(p)}}}{\|\mathbf{B^{(p)}}\|}$.

Alternatively, the optimal $\mathbf{A_p^*}$ can be equivalently
obtained by calculating the eigenvector corresponding to the largest
eigenvalues of the following matrix:
\begin{equation}\label{4.4.7}
\sum_{i=1}^{NM}\sum_{j=1}^{NM}exp\left(\mu_i+\mu_j+\sum_{k=1}^{p-1}\left(C_{ik}^*+C_{kj}^*\right)
\right)\mathbf{C^{Ch}_{i\cdot}}\mathbf{C^{Ch}_{\cdot
j}}=diag(\mathbf{d^{(p)}})\Sigma_{MN}diag(\mathbf{d^{(p)}})
\end{equation}
\noindent and after imposing the orthonormality condition by the QR
method.

%
%

\section{Simulation Framework}

We consider  the constant volatility case only, and run our
simulation with different combinations of path-generation techniques
and different random number generators.

As far as path-generation methods are concerned we use the standard
Cholesky, the PCA and the two LT decompositions for Asian options
introduced in the previous subsections. In particular for the first
two approaches we rely on the properties of the Kronecker product in
order to compute the decomposition fast (see Dahl,  Benth
\cite{DB2001} and \cite{DB2002} and Sabino \cite{Sab2007} for
further details).

LT methods require the iterative calculations of orthogonal
matrices. We attain the task implementing an \emph{ad hoc} QR
factorization, as described in Appendix, that does not require high
computational cost. For the LT decomposition the total computational
time is than the sum of the time to compute the Cholesky and the
optimal orthogonal matrix $A$.

The numerical test consists of three main steps:
\begin{enumerate}
  \item Random number generation by standard MC, LHS or RQMC.
  \item Path generation with Cholesky, PCA, and the two LT algorithm discussed above (LT1 and LT2, respectively).
  \item MC estimation.
\end{enumerate}

As RQMC generator we use a Faure-Tezuka scrambled version of the
$50$-dimensional Sobol´ sequence satisfying Sobol´s property A (see
Glasserman \cite{Glass2004}, J\"{a}ckel \cite{Ja2003} and Owen
\cite{ow2002} for further details). We pad the remaining random
components out with LHS. This strategy is intended to investigate
the effective improvement of the LT methods when coupled with QMC.
Indeed, it can be proven that the LHS gives good variance reductions
when the target function is sum of one-dimensional functions (see
Glasserman \cite{Glass2004} and Owen \cite{ow1998B}). On the other
hand, the LT methods is conceived to capture the lower effective
dimension in superposition sense for linear combinations. As a
consequence, we should already observe a high accuracy when running
the simulation with LHS and LT. Our setting is thought to test how
large is the improvement given by the LT factorization. We compute a
suboptimal $A$ up to dimension $50$ in order to be coherent with the
choice of the $50$-dimensional Sobol´ sequence.

Stratification introduces correlation among random drawings so
that the hypothesis of the Central Limit Theorem are not satisfied
and we cannot compute the RMSE straightforward.
We rely on the batch methods that consists of repeating $%
N_{B}$ simulations for $B$ times (batches).
%
%
\section{Numerical Investigations}

\noindent We develop our simulation procedure in order to test the
computational burden and the efficiency of the Linear Transform
method. We compare its results with those obtained with standard
techniques like Cholesky and PCA decompositions. Furthermore, we use
several random number generators, in particular, we adopt a
Faure-Tezuka scrambled version of the $50$-dimensional Sobol'
sequence satisfying the Sobol's property A.

As a numerical example, we estimate the fair price of an Asian
option on a basket of $M=10$ underlying assets with $N=250$
sampled points.

The chosen parameters are those in the original paper of Imai and
Tan \cite{IT2007} and are shown in Table \ref{test_S2}.
\begin{table}[tbp] \centering%
\begin{tabular}{||c||}
\hline\hline
$%
\begin{array}{lll}
S_{i}\left( 0\right) & = & 100 \\
K & = & 90,100 \text{ and } 110 \\
r & = & 4\% \\
T & = & 1 \\
\sigma_i & = & 10\%+\frac{i-1}{9}40\%\quad \text{for }i=1,\dots,10
\\\rho_{ij} & = & 0\quad \textrm{and}\quad40\%\quad \text{for }i,j=1,\dots,10
\end{array}%
$ \\ \hline\hline
\end{tabular}%
\caption{Input Parameters.}\label{test_S2}%
\end{table}

The nominal dimension of the problem is $M\times N=2500$
equal to the number of rows and columns of the global correlation matrix $%
\Sigma_{MN}$.

We perform the path-generation by computing the Cholesky, the PCA
and two versions of the LT decompositions of the global correlation
matrix $\Sigma_{MN}$. We label LT2 for the general case and LT1 for
the method described for Asian options only. As far as the first two
generations are concerned, we rely on the properties of the
Kronecker product in order to reduce the computational burden as
described in Dahl, Benth \cite{DB2001} and \cite{DB2002} and Sabino
\cite{Sab2007}.

As far as the implementation of the two LT methods proposed by Imai
and Tan is concerned, we apply the fast version of the QR
decomposition described in the appendix.

The simulation procedure is implemented in MATLAB running on a
laptop with an Intel Pentium M, processor 1.60 GHz and 1 GB RAM.

Table \ref{Var_Contr} shows the percentage of the cumulative
contribution of the variance for the first $10$ components both for
the zero and positive correlation cases, this is the ratio between
equation (\ref{4.4.2}) and (\ref{4.4.1}) with $p$ up to $10$.
\begin{table}[tbp] \centering%
\begin{tabular}{||c|cccc|cccc||}
\hline\hline
 &  & Uncorrelation & &  &  & Correlation &  & \\
 \hline
Dimension & Cholesky & PCA & LT1 & LT2 & Cholesky & PCA & LT1 & LT2\\
\hline
1 & 0.01 & 23.15 & 88.41 & 88.41 & 0.41 & 91.40 & 94.41 & 94.41\\
2 & 0.03 & 42.24 & 90.68 & 90.24 & 0.60 & 93.04 & 95.27 & 95.17\\
3 & 0.07 & 57.69 & 93.53 & 92.20 & 0.71 & 94.32 & 96.25 & 95.91\\
4 & 0.13 & 69.91 & 95.10 & 94.11 & 0.79 & 95.28 & 97.02 & 96.58\\
5 & 0.21 & 79.30 & 96.69 & 95.45 & 0.85 & 95.99 & 97.59 & 97.25\\
6 & 0.31 & 86.25 & 97.83 & 96.68 & 0.91 & 97.42 & 97.93 & 97.58\\
7 & 0.45 & 91.14 & 98.32 & 97.53 & 0.95 & 97.93 & 98.34 & 98.04\\
8 & 0.61 & 94.33 & 98.60 & 98.10 & 1.00 & 98.28 & 98.61 & 98.24\\
9 & 0.81 & 94.87 & 98.85 & 98.45 & 1.04 & 98.51 & 98.75 & 98.45\\
10 & 1.05 & 95.28 & 98.97 & 98.64 & 1.08 & 98.78 & 98.92 & 98.61\\
\hline\hline
\end{tabular}%
\caption{Percentage of Variance Contribution up to dimension $10$.}\label{Var_Contr}%
\end{table}
All results up to $p=5$ are consistent with those presented by Imai
and Tan \cite{IT2007}. It can be noticed that the LT is the best
performing path-generation technique in the statistical sense
specified above, where the first specification is a bit better. The
PCA decomposition is  almost as accurate as the LT approach for the
correlation case only.

The effective dimensions found with each method are reported in
Table \ref{Eff_Dim}. The advantage of the PCA and LT methods with
respect to the Cholesky decomposition is evident both for the
correlation and uncorrelation cases. The Cholesky decomposition
collects $98.58\%$ and $98.70\%$ of the total variance for  $p=2000$
for the uncorrelation and correlation cases, respectively.
\begin{table}[tbp] \centering%
\begin{tabular}{||cccc|cccc||}
\hline\hline
  & Uncorrelation & &  &  & Correlation &  & \\
 \hline
 Cholesky & PCA & LT1 & LT2 & Cholesky & PCA & LT1 & LT2\\
\hline
 $d_T > 2000$  & $d_T = 18$ & $d_T = 11$ & $d_T = 13$ & $d_T > 2000$ & $d_T = 12$ & $d_T = 12$ & $d_T = 12$\\
\hline\hline
\end{tabular}%
\caption{Effective Dimensions.}\label{Eff_Dim}%
\end{table}

We compute the computational times elapsed to decompose the global
covariance matrix with each method so that we can compare the
efficiency of all the methods; we compute only $50$ optimal
columns for the LT technique. Table \ref{El_Tim} shows the
estimated times in seconds.

\begin{table}[tbp] \centering%
\begin{tabular}{||c|cccc|cccc||}
\hline\hline
 & & Uncorrelation & &  &  & Correlation &  & \\
 \hline
 & Cholesky & PCA & LT1 & LT2 & Cholesky & PCA & LT1 & LT2\\
\hline
 $time$ & $0.60$  & $25.77$ & $53.14$ & $53.21$ & $0.59$ & $25.55$ & $53.02$ & $53.20$\\
\hline\hline
\end{tabular}%
\caption{Computational Times.}\label{El_Tim}%
\end{table}
The computational times we found are a lot lower than those
presented by Imai and Tan \cite{IT2007} despite the fact that are
computed with a slower computer. In particular, the implementation
of the LT method with the QR approach presented in the appendix (up
to $50$ columns) is more efficient of a factor tirthy. Furthermore,
the LT methods has the versatility to allow the computation of a
suboptimal matrix that is statistically justified by ANOVA
considerations. In contrast, the PCA decomposition lacks this
possibility without losing information.

In the case of time-depending volatilities we could not rely on the
properties of the Kronecker products in order to reduce the
computational costs to run the PCA decomposition of the global
covariance matrix. In contrast, the \emph{ad hoc} QR approach for
the LT method is preserved and needs a computational time of the
same order as we will present in future studies.

In the case of time-dependent volatilities it is fundamental to
implement a fast Cholesky decomposition to be coupled with the QR
method (see Sabino \cite{Sab2007} for further details of this type
of Cholesky algorithm).

As a final step, we launch a MC simulation in order to estimate the
fair price of the Asian basket option with $8192$ generations and
$10$ replications.

As already mentioned, we use a standard pseudo-random generator,
the LHS method and a $50$-dimensional Faure-Tezuka scrambled
version of the Sobol' sequence satisfying the Sobol's property A.

It is known that (R)QMC simulations do not yield any improvements
with respect to standard MC ones when the problem dimension is high
(generally d $\ge 20/30$). Owen \cite{ow1998B} proposes mainly two
approaches to extend the better convergence of the (R)QMC in high
dimensions: the Latin Super Cube method and the padding with LHS.

Briefly, the former  consists of grouping the input variables and
rearranging their order  with a random permutation. The latter
consists in fixing the more important variables and then pad the
remaining ones out with the LHS.

Even if this last method requires more computational costs, it can
give further insight into the LT method. Indeed, it can test if the
LT really selects the best variables in statistical sense and
reduces the effective dimension. For the presented case we compare
its results with those obtained with the pure LHS generator. We
then, choose a $50$-dimensional Sobol' sequence, coherent with the
suboptimal matrix $A$, and pad the remaining $2450$ dimensions out.
\begin{table}[tbp] \centering
\begin{tabular}{||c|cc|cc|cc||}
\hline\hline
&&&Standard MC&&&\\
\hline
&K=90&&K=100&&K=110&\\
\hline
&Price&RMSE&Price&RMSE&Price&RMSE\\
\hline
Cholesky&$12.257$&$0.038$&$5.604$&$0.029$&$2.007$&$0.018$\\
PCA&$12.291$&$0.038$&$5.648$&$0.029$&$2.040$&$0.018$\\
LT1=50&$12.240$&$0.038$&$5.681$&$0.029$&$2.004$&$0.018$\\
LT2=50&$12.256$&$0.038$&$5.687$&$0.029$&$2.007$&$0.018$\\
\hline
&&&LHS&&&\\
\hline
&Price&RMSE&Price&RMSE&Price&RMSE\\
\hline
Cholesky&$12.3320$&$0.0097$&$5.670$&$0.013$&$2.0394$&$0.0084$\\
PCA&$12.3292$&$0.0025$&$5.6655$&$0.0032$&$2.0389$&$0.0034$\\
LT1=50&$12.3291$&$0.0015$&$5.5968$&$0.0019$&$2.0332$&$0.0022$\\
LT2=50&$12.2468$&$0.0022$&$5.6015$&$0.0017$&$2.0348$&$0.0017$\\\hline
&&&RQMC&&&\\
\hline
&Price&RMSE&Price&RMSE&Price&RMSE\\
\hline
Cholesky&$12.3410$&$0.0094$&$5.631$&$0.014$&$2.021$&$0.012$\\
PCA&$12.32900$&$0.00060$&$5.65770$&$0.00039$&$2.03360$&$0.00041$\\
LT1=50&$12.32800$&$0.00036$&$5.65720$&$0.00040$&$2.03400$&$0.00021$\\
LT2=50&$12.32800$&$0.00025$&$5.65690$&$0.00019$&$2.03420$&$0.00039$\\
\hline\hline
\end{tabular}
\caption{Correlation Case: Estimated Prices and Errors.}\label{Corr_S}%
\end{table}
\begin{table}[tbp] \centering
\begin{tabular}{||c|cc|cc|cc||}
\hline\hline
&&&Standard MC&&&\\
\hline
&K=90&&K=100&&K=110&\\
\hline
&Price&RMSE&Price&RMSE&Price&RMSE\\
\hline
Cholesky&$11.560$&$0.021$&$3.426$&$0.015$&$0.3625$&$0.0051$\\
PCA&$11.553$&$0.021$&$3.414$&$0.015$&$0.3591$&$0.0052$\\
LT1=50&$11.454$&$0.021$&$3.356$&$0.015$&$0.3525$&$0.0051$\\
LT2=50&$11.454$&$0.021$&$3.357$&$0.015$&$0.3523$&$0.0051$\\
\hline
&&&LHS&&&\\
\hline
&Price&RMSE&Price&RMSE&Price&RMSE\\
\hline
Cholesky&$11.5915$&$0.0037$&$3.432$&$0.007$&$0.3605$&$0.0038$\\
PCA&$11.5913$&$0.0064$&$3.4546$&$0.0054$&$0.3686$&$0.0030$\\
LT1=50&$11.4754$&$0.0013$&$3.3666$&$0.0023$&$0.3662$&$0.0021$\\
LT2=50&$11.4779$&$0.0030$&$3.3693$&$0.0022$&$0.3662$&$0.0015$\\
\hline
&&&RQMC&&&\\
\hline
&Price&RMSE&Price&RMSE&Price&RMSE\\
\hline
Cholesky&$11.5890$&$0.0033$&$3.4351$&$0.0049$&$0.3605$&$0.0035$\\
PCA&$11.59000$&$0.00039$&$3.4444$&$0.0015$&$0.3662$&$0.0011$\\
LT1=50&$11.59100$&$0.00025$&$3.44360$&$0.00039$&$0.36673$&$0.00034$\\
LT2=50&$11.59200$&$0.00036$&$3.44440$&$0.00033$&$0.36599$&$0.00034$\\
\hline\hline
\end{tabular}
\caption{Uncorrelation Case: Estimated Prices and Errors.}\label{UnCorr_S}%
\end{table}

Tables \ref{Corr_S} and \ref{UnCorr_S} show the results of our
numerical experiment. All values are  statistically consistent, but
exhibit a different accuracy.

As expected, standard Cholesky decomposition is almost not sensitive
to the used random generation technique and gives the worst results.

LT and PCA decompositions provide good improvements for the RMSEs
for both the RQMC and LHS generations. As far as the last method is
concerned, we note that it is sensitive to the decomposition used
and returns lower RMSEs when the LT decomposition is applied. This
means that the LT approach is really reducing the effective
dimension in superposition sense, "splitting" the integrand function
into a sum of linear functions.

As already mentioned, the LHS should reduce the RMSE in the case the
integrand function is the sum of one-dimensional functions. This is
best accomplished by the LT as evident from the above results.

The RQMC simulation and the LT decompositions confirm their superior
performance.

It can be noted that the RQMC is sensitive to the used decomposition
approach and does not have any advantage over the LHS when we use
the Cholesky decomposition.

Our evaluations return RMSEs with the same accuracy as Imai and Tan
\cite{IT2007} when we only consider a $50$-dimensional Sobol'
sequence without using the complete LSS.

Our framework is more extreme and the LT provides the same
efficiency for all the strike prices and all correlations
considered. In contrast, the PCA approach gives high improvements
only in the correlation case.

The general  and the Asian options settings of the LT decompositions
are almost equally performing with the latter one giving slightly
better results.

We can conclude that the LT is the best decomposition method and
tremendously enhances QMC simulations because it optimally reduces
the effective dimension of the problem.

The LT construction can be made faster from the computational point
of view, provided we implement the QR decomposition described in the
appendix.
\section{Conclusion}
In this paper we investigate the accuracy of the LT, introduced by
Imai and Tan, both from the computational and the accuracy points of
view. In particular, we implement a numerical procedure based on the
QR factorization that realize the LT decomposition fast. Moreover,
we extensively investigate the improvements the LT gives to QMC
methods that is sensitive to the effective dimension of the problem.

As a numerical test we launch a high-dimensional simulation with the
same set of parameters as in Imai and Tan \cite{IT2007} in order to
price Asian basket options.

Our setting is more extreme than the one discussed in the cited
references. We do not rely on the complete LSS high-dimensional
extension of the features of the QMC but we use a lower dimensional
scrambled Sobol' sequence only, and pad the remaining ones out with
LHS.

We compare these results with those published by Imai and Tan
\cite{IT2007} and those we found when using different decompositions
and different random number generators.

The LT construction provides the best accuracy with respect to the
standard Cholesky approach and the PCA decomposition.

It provides considerable improvements even when simulations are
carried out with a partial RQMC method. The LT accuracy is still
notably better than the one we found with the complete LHS. In
particular, we attain RMSEs of the same order as those presented by
Imai and Tan.

Moreover, the fast QR decomposition we implement gives an
improvement of a factor $30$ in terms of computational time compared
to the results presented in Imai and Tan \cite{IT2007} calculated
with a slower computer.

PCA decomposition enhances QMC simulations but still requires a high
computational burden when time-dependent volatilities are considered
and does not give the versatility to find a suboptimal matrix
without introducing bias (see Sabino \cite{Sab2007} for details) .

Our QR-implementation makes the LT more efficient and
computationally more convenient while maintaining its versatility
for different problems.

\newpage
\section{Appendix}
\subsection{The QR Method}
The QR factorization of an $m$-by-$n$ matrix $A$ is given by:
\begin{equation}\label{A.1}
A = QR
\end{equation}
\noindent where $Q\in \mathbb{R}^{m\times m}$ is orthogonal and
$R\in \mathbb{R}^{m\times n}$ is upper triangular. A fundamental
result is that if $A$ has full column rank, then the first $n$
columns of $Q$ form an orthonormal basis of $ran(A)$. As a
consequence, the QR factorization provides a way to return an
orthonormal basis for a set of (independent) vectors. Different
approaches can be chosen to calculate the QR decomposition such as
the Householder and Givens transformations (see Golub, Van Loan
\cite{GV1996} as a fundamental reference).

The former transformations are rank-two corrections of the
identity of the form:
\begin{eqnarray}\label{A.G.}
    G(i,k,\theta)&=&\left[
  \begin{array}{ccccccc}
  1 & \cdots & 0 & \cdots & 0 & \cdots & 0\\
  \vdots & \ddots & 0 & \vdots & 0 & \vdots & 0\\
  0 & \cdots & \cos\theta & \cdots & \sin\theta & \dots & 0\\
  \vdots & \cdots & 0 & \cdots & 0 & \dots & 0\\
  0 & \cdots & -\sin\theta & \cdots & \cos\theta & \dots & 0\\
  \vdots & \cdots & 0 & \cdots & 0 & \dots & 0\\
  0 & \cdots & 0 & \cdots & 0 & \dots & 1\\
\end{array}\right]
\begin{array}{c}
\\
\\
i\\
\\
k\\
\\
\\
\end{array}
\\
 & &
\begin{array}{cccccccccccc}
    &&&&&&i&&&&&k
\end{array}
\end{eqnarray}
$G(i,k,\theta)$ performs a counterclockwise rotation of $\theta$
radians in the $(i,k)$ plane.

Householder and Givens transformations are orthogonal
transformations constructed in order to introduce zeros in a vector.
Indeed, suppose we are given with $\mathbf{0}\neq\mathbf{x}\in
\mathbb{R}^n$ we can find $H$ and $G$, the Householder and Given
Rotation, respectively, that annihilate the $k$-th component of
$\mathbf{x}$:
\begin{eqnarray}
  H^T \mathbf{x} &=& \mathbf{y},y_k=0 \\
  G^T \mathbf{x} &=& \mathbf{y},y_k=0
\end{eqnarray}
The following scheme illustrates the idea for QR factorization with
Givens rotations:
\begin{eqnarray*}A=\left[
  \begin{tabular}{ccc}
  x & x & x \\
  0 & x & x \\
  0 & x & x \\
  0 & 0 & x \\
\end{tabular}\right] \underrightarrow{(2,3)}
G_1^TA=\left[\begin{tabular}{ccc}
  x & x & x \\
  0 & x & x \\
  0 & 0 & x \\
  0 & 0 & x \\
\end{tabular}\right] \underrightarrow{(3,4)}
G_2^TG_1^TA=R
\end{eqnarray*}

\noindent Here we have highlighted the 2-vectors that define the
underlying Given rotations. Generally $Q=G_1\cdots G_n$ where $n$
is the total number of rotations and $R=Q^TA$.

Consider $A\in\mathbb{R}^{m\times n}$, its QR decomposition
$A=Q_AR_A$ and $B\in\mathbb{R}^{m\times (n+1)}$ with the first $n$
columns equal to $A$, $Q_B$ and $R_B$. The QR decomposition of $B$
can be easily obtained from $Q_A$ and $R_A$.

Denote $\mathbf{b}$ the last column of $B$, so that $B = [A $
$\mathbf{b}]$ it leads to  $Q_A^T B = [R_A$ $ Q_A^T\mathbf{b}] =
\tilde{R}_B$. $\tilde{R}_B$ has the following form:
\begin{equation}
 \tilde{R}_B =
    \left[\begin{array}{cccc}
            x & x & x & x \\
            0 & x & x & x \\
            0 & 0 & x & x \\
            0 & 0 & x & x \\
            0 & 0 & x & x
          \end{array}
       \right].
\end{equation}
\noindent In order to obtain the complete QR factorization of $B$ we
only need to find $t$ Givens transformations $G_1,\dots,G_t$ that
introduce zeros in the $n$-th column making $R_B
=G_t^T,\dots,G_1^T\tilde{R}_B$ upper triangular.

Summarizing $Q_B =G_1,\dots,G_tQ_A$ and $R_B
=G_t^T,\dots,G_1^TQ_AB$.


\end{document}